\documentclass[leqno,12pt]{amsart}
\usepackage{amsmath,amssymb,amsthm,amsfonts,amscd}
\usepackage[all]{xy}
\usepackage[dvips]{graphicx}

\newcommand{\imm}{\looparrowright}

\newcommand{\R}{\mathbb{R}}
\newcommand{\Z}{\mathbb{Z}}
\newcommand{\N}{\mathbb{N}}

\newcommand{\co}{\colon\thinspace}

\newtheorem{thm}{Theorem}[section]
\newtheorem{thm*}{Theorem}

\newtheorem{prop}[thm]{Proposition}
\newtheorem{defn}[thm]{Definition}

\newtheorem{lemma}[thm]{Lemma}

\newtheorem*{rems}{Remarks}

\begin{document}

\title[Self-intersections]{Self-intersections of immersions\\
and Steenrod operations}

\author[P.J. ~Eccles, M. ~Grant]{Peter J. Eccles and Mark Grant}
\address{School of Mathematics,
    Alan Turing Building, The University of Manchester, Oxford Road, Manchester, M13 9PL, UK}
\address{School of Mathematical Sciences, The University of Nottingham,
University Park, Nottingham, NG7 2RD, UK}
\email{peter.eccles@manchester.ac.uk}
\email{mark.grant@nottingham.ac.uk}

\maketitle
\begin{abstract}
We present a formula describing the action of a generalised Steenrod operation of $\Z_2$-type \cite{tD68} on the cohomology class represented by a proper self-transverse immersion $f\co M\imm X$. Our formula depends only on the Umkehr map, the characteristic classes of the normal bundle, and the class represented by the double point immersion of $f$. This generalises a classical result of R.\ Thom \cite{Th50}: If $\alpha\in H^k(X;\Z_2)$ is the ordinary cohomology class represented by $f\co M\imm X$, then $Sq^i(\alpha)=f_* w_i(\nu_f)$.
\end{abstract}

\section{Introduction}

Let $f\co M\imm X$ be a proper self-transverse immersion (all manifolds and maps are assumed to be differentiable of class $C^{\infty}$). Then the self-intersection set of $f$,
\[
\Sigma = \{(x,y)\in M\times M\mid x\neq y, f(x)=f(y)\},
\]
is a manifold, and $f$ induces a proper immersion $\psi\co \Sigma\imm X$ with image the double points of $f$. The cyclic group $G=\Z_2$ acts on $\Sigma$ by swapping factors and trivially on $X$, and $\psi$ is $G$-equivariant. Passing to orbit spaces we obtain a proper immersion $\Psi\co\Sigma/G \imm X$, the so-called {\em double point immersion of $f$}.

On the other hand, let $h^*$ be a multiplicative generalised cohomology theory for which the normal bundle $\nu_f$ of $f\co M\imm X$ is $h^*$-oriented. Then $f$ represents a class $f_*(1) \in h^*(X)$. Let $E\to B$ be a principal $G$-bundle in which both $E$ and $B$ are connected manifolds. We assume that the theory $h^*$ admits an internal Steenrod operation $\mathcal{P}\co h^{d*}(-)\to h^{2d*}(B\times -)$ of type $(G,d)$ as defined by T.\ tom Dieck \cite{tD68} (more details will be given in Section 4 below), and that $d$ divides the codimension of $f$. Then the normal bundle $\nu_{\Psi}$ inherits an $h^*$-orientation from that of $\nu_f$, and hence the double point immersion represents a class $\Psi_*(1)\in h^*(X)$. Our main result may now be stated as follows.

\begin{thm*}\label{main}
Let $f\co M^{n-dk}\imm X^n$ be a proper self-transverse immersion representing $f_*(1)\in h^{dk}(X)$, and assume that $h^*$ admits a Steenrod operation $\mathcal{P}$ of type $(G,d)$. Then
\[
\mathcal{P}\big(f_*(1)\big) = 2\times\Psi_*(1) + (f\times 1)_*e(\gamma\otimes\nu_f) \in h^{2dk}(B\times X).
\]
Here $e$ is the Euler class, $\gamma$ is the line bundle associated to $E\to B$, and $\otimes$ denotes exterior tensor product of vector bundles.
\end{thm*}

This formula can be viewed as an equivariant version of Whitney's formula $f_*(1)\cup f_*(1) = 2\Psi_*(1) + f_*e(\nu_f)$, obtained by pushing forward into $X$ the double point case of Herbert's formula \cite{EG07}. It should be compared with the results in Section 3 of the landmark paper of D.\ Quillen \cite{Qu71}, where a formula is derived relating the two types of operation (Steenrod and Landweber-Novikov) in complex cobordism. In a certain sense, our formula shows that the difference between these two types of operation is precisely measured by the cobordism class of the double point immersion, at least on classes represented by immersions. We also mention Theorem 7.1 in the book of Buoncriatiano, Rourke and Sanderson \cite{BRS}, which is similar in spirit.

The proof employs sub-cartesian diagrams, as introduced by Felice Ronga \cite{Ron} (and subsequently used in \cite{BL}). These generalise from embeddings to immersions the clean intersection formula used by Quillen in \cite{Qu71}, removing the need for stabilisation and hence leading to sharper results in the immersed case. After some preliminaries on orientations and Umkehr maps in Section 2, we recall the necessary facts about sub-cartesian diagrams in Section 3. Section 4 reviews material on Steenrod operations in generalised cohomology theories, most of which can be found in \cite{tD68}. In Section 5 we prove Theorem \ref{main}. Finally, in Section 6 we observe that in ordinary mod $2$ cohomology our formula reduces to a classical result normally attributed to R.\ Thom \cite{Th50}, namely that the action of the Steenrod squares on a class represented by an immersion is given in terms of the Stiefel-Whitney classes of the normal bundle by the formula $Sq^i(f_*(1))=f_* w_i(\nu_f)$.

Theorem \ref{main} corrects an erroneous formula in the second author's thesis (see \cite[Theorem 6.2]{G} and its corollaries). We are grateful to Sergey Melikhov for pointing out the mistake (which occurs in the proof of \cite[Proposition 5.9]{G}), and for his continued interest in this work.

\section{Orientations and Umkehr maps}

We first recall some facts about orientations and Umkehr maps in generalised cohomology theories, for which the book by E.\ Dyer \cite{D} is a standard reference. Let $h^*$ be a multiplicative generalised cohomology theory on the category of topological pairs $(X,A)$. For any space $X$ we denote $h^*(X)=h^*(X,\emptyset)$, which is a graded ring with unit $1\in h^0(X)$. If $\xi$ is a $k$-dimensional vector bundle over a paracompact space $X$, its {\em Thom space} is the pointed space $T\xi:= D\xi/S\xi$, where $D\xi$ and $S\xi$ are the disk and sphere bundles of $\xi$ with respect to some Riemannian metric. An {\em $h^*$-orientation} of $\xi$ is a choice of Thom class $t_\xi\in h^k(T\xi,*)\cong h^k(D\xi,S\xi)$. The {\em Euler class} of an $h^*$-oriented vector bundle $\xi$ is the class $e(\xi)=i^*(t_\xi)\in h^k(X)$, where $i\co X\to D\xi$ is the zero section. If $f\co X'\to X$ is a map, the pullback bundle $f^*\xi$ over $X'$ receives an induced orientation from that of $\xi$.
\begin{lemma}
 Let $\xi$ and $\eta$ be bundles over $X$. If any two of the bundles $\xi$, $\eta$ and $\xi\oplus\eta$ are $h^*$-oriented, then so is the third.
\end{lemma}
\begin{proof}
See \cite[Chapter 1, Section C]{D}.
\end{proof}

Let $f\co M^{n-k}\imm X^n$ be a proper immersion. Assuming that $X$ has been given a metric, the {\em normal bundle} $\nu_f$ of $f$ is defined by $f^*TX\cong TM\oplus \nu_f$. We say $f$ is {\em $h^*$-oriented} if its normal bundle $\nu_f$ is $h^*$-oriented. Then $f$ and its orientation induce an {\em Umkehr} (or {\em Gysin} or {\em pushforward}) map
\[
f_*\co h^*(M)\to h^{*+k}(X).
\]
We refer the reader to \cite[Chapter 1, Section D]{D} for more details.
\begin{lemma}\label{Um}
The Umkehr maps satisfy the following properties:
\begin{enumerate}
\item If $g\co N\imm M$ and $f\co M\imm X$ are proper $h^*$-oriented immersions, and $\nu_{f\circ g}\cong \nu_g\oplus g^*\nu_f$ is given the induced $h^*$-orientation, then
    \[
    (f\circ g)_* = f_*\circ g_*.
    \]
\item Let $f\co M\imm X$ be a proper $h^*$-oriented immersion, and let $N$ be a manifold. For any $\alpha\in h^*(N)$ and $\beta\in h^*(M)$ we have
    \[
    (1\times f)_*(\alpha\times\beta) = \alpha\times f_*(\beta)\in h^*(N\times X).
    \]
\item If $p\co \tilde{X}\imm X$ is an $n$-sheeted covering of manifolds, then the composition $p_*\circ p^*\co h^*(X)\to h^*(X)$ is multiplication by $n$.
\end{enumerate}
\end{lemma}
\begin{proof}
See Dyer \cite[Chapter 1, Section D]{D} or J.\ F.\ Adams \cite[Chapter 4]{A}. Alternatively, the proofs can be deduced quite easily from Proposition \ref{Ronga} in the next Section.
\end{proof}
\begin{defn}
If $f\co M^{n-k}\imm X^n$ is a proper $h^*$-oriented immersion, we say that $f$ {\em represents} the class $f_*(1)\in h^k(X)$, where $1\in h^0(M)$ is the unit.
\end{defn}

\section{Sub-cartesian diagrams}

In this Section we recall the definition of sub-cartesian diagrams and
a Proposition due to F.\ Ronga \cite{Ron}. This is a generalisation of the notion of clean intersection of submanifolds, defined by R.\ Bott \cite{Bo} and used by Quillen \cite{Qu71}.

\begin{defn}\label{clean}
The commutative diagram of proper immersions
\[
 \xymatrix{
 Z \ar[r]^-{\beta} \ar[d]_{\alpha}  & B \ar[d]^{f} \\
 A \ar[r]^-{g}         &  X}
 \]
is called {\em sub-cartesian} if
 \begin{itemize}
 \item the map $(\alpha,\beta)\co Z \to A\times B$ is an embedding onto $$A\times_X B = \{ (a,b)\in A\times B \mid g(a) = f(b)\};$$
 \item the following sequence of bundles over $Z$ is exact,
 \[
 \xymatrix{
0 \ar[r] &  TZ \ar[r]^-{d\alpha\oplus d\beta} & \alpha^*TA\oplus \beta^*TB \ar[r]^-{dg -df} & \beta^*f^*TX, }
\]
where $(dg - df)(v,w) = dg(v) - df(w)$.
\end{itemize}
The bundle $E=\mathrm{coker}(df - dg)$ is called the {\em excess bundle}.
\end{defn}
\begin{rems}{\em
The first condition says that the manifold $Z$ is the intersection of $f$ and $g$, where multiple points are counted with the relevant multiplicity. The second condition describes the tangent space of $Z$ locally as the intersection of the tangent spaces of $A$ and $B$ in $X$. That is to say, $f$ and $g$ intersect cleanly, in the terminology of Quillen \cite{Qu71}.

 Note that $Z$ is not assumed to be of constant dimension, and so neither is $E$ in general.

 The excess bundle $E$ is the zero bundle if and only if $f$ and $g$ are transverse maps.}
\end{rems}

The second condition in Definition \ref{clean} implies that there is an exact sequence of bundles
\[
 \xymatrix{
  0 \ar[r] & \nu_\alpha \ar[r] & \beta^*\nu_f \ar[r] & E \ar[r] & 0}
 \]
 over $Z$, and hence (after choosing metrics) an isomorphism of bundles
 \[
 \beta^*\nu_f  \cong \nu_\alpha \oplus E.
 \]
\begin{prop}[F.\ Ronga \cite{Ron}]\label{Ronga} Given a sub-cartesian diagram
\[
 \xymatrix{
 Z \ar[r]^-{\beta} \ar[d]_{\alpha}  & B \ar[d]^{f} \\
 A \ar[r]^-{g}         &  X}
 \]
suppose that $\nu_f$ and one of the bundles $E$ or $\nu_\alpha$ are $h^*$-oriented. Then for any $c\in h^*(B)$ we have
\begin{equation}
g^*f_*(c)=\alpha_*(e(E)\cdot\beta^*(c))\in h^*(A),
\end{equation}
where the third orientation is determined by $\beta^*\nu_f  \cong \nu_\alpha \oplus E$.
\end{prop}

\section{Steenrod operations}

Steenrod operations arise from higher homotopy
commutativity of the product in a multiplicative generalised
cohomology theory. They were originally defined in $\Z_p$-cohomology by
Steenrod \cite{SE}, and later in K-theory by
Atiyah \cite{At66}. The so called Steenrod-tom Dieck cobordism operations, defined in \cite{tD68}, are a key tool in Quillen's work relating complex cobordism and the theory of formal groups. Here we recall some details of the axiomatic treatment of Steenrod operations given in \cite{tD68}.

We restrict our attention to $G=\Z_2$, the cyclic group of order two, although the statements in this Section admit analogues for any subgroup of a symmetric group. We fix a principal $G$-bundle $E\to B$ in which both total and base space are connected manifolds. Given any $G$-space $Y$ we may form the Borel construction $Y_G:=E\times_G Y$, which is a manifold if $Y$ is a $G$-manifold. This construction is functorial with respect to $G$-maps. In particular, to a given a $G$-vector bundle $\xi\to Y$ we may associate a vector bundle $\xi_G\to Y_G$ of the same dimension. Note that choosing a basepoint $e\in E$ results in a map $i\co Y\to Y_G$, and $i^*\xi_G\cong\xi$.

As a particular case of these constructions, an arbitrary space $X$ can be regarded as a $G$-space with trivial action, whilst $G$ acts on the product $X\times X$ by transposition of factors. The diagonal map $\triangle(x)=(x,x)$ is evidently equivariant, and we have a diagram
\[
\xymatrix{
  & X \ar[r]^{\triangle} \ar[ld]_{\iota} \ar[d]_{i} & X\times X \ar[d]_{i} \\
 B\times X \ar@{=}[r]&  X_G \ar[r]^-{\triangle_G} & (X\times X)_G}
\]
We are now ready to define the Steenrod operations.

\begin{defn}
Let $d\geq 1$ be a natural number. An {\em external Steenrod operation} of type $(G,d)$ in the multiplicative cohomology theory $h^*$ is a sequence $P$ of natural transformations
\[
P^{dk}\co h^{dk}(X)\to h^{2dk}((X\times X)_G),\quad k\in \Z
\]
such that $i^* P^{dk}(\alpha) = \alpha\times\alpha \in h^{2dk}(X\times X)$.

An {\em internal Steenrod operation} of type $(G,d)$ in $h^*$ is a sequence $\mathcal{P}$ of natural transformations
\[
\mathcal{P}^{dk}\co h^{dk}(X)\to h^{2dk}(B\times X),\quad k\in \Z
\]
such that $\iota^* \mathcal{P}^{dk}(\alpha) = \alpha\cup\alpha\in h^{2dk}(X)$.
\end{defn}

Note that an internal Steenrod operation arises from an external one by setting $\mathcal{P}=\triangle_G^*\circ P$. In \cite{tD68}, tom Dieck constructs Steenrod operations of type $(G,d)$ in each of the cobordism theories $h^*=M\Gamma^*$, where $\Gamma=O,SO,U,SU,Sp$ and $d=1,2,2,4,4$ respectively.

Suppose $h^*$ admits a Steenrod operation $P$ of type $(G,d)$. Let $f\co M\imm X$ be a proper $h^*$-oriented immersion representing $f_*(1)\in h^{dk}(X)$. Then the class $P\big(f_*(1)\big)$ is also represented by an immersion, as follows (this geometric construction of $P$ is due to tom Dieck \cite{tD68} and Quillen \cite{Qu71}). Given any $h^*$-oriented bundle $\xi\to X$, the product bundle $\xi\times\xi\to X\times X$ may be given the product orientation. This is also a $G$-bundle under transposition of factors, and the associated bundle
\[
(\xi\times \xi)_G\to (X\times X)_G
\]
has a canonical $h^*$-orientation (a Thom class is given by $\tilde{P}t_\xi$, where $\tilde{P}$ is a reduced Steenrod operation; see \cite[Satz 4.1]{tD68}). In particular, the bundle
\[
(\nu_f\times\nu_f)_G \to (M\times M)_G
\]
is canonically $h^*$-oriented. This is the normal bundle of the immersion $(f\times f)_G\co (M\times M)_G\imm (X\times X)_G$, which is thus $h^*$-oriented.

\begin{prop}[\cite{tD68}, \cite{Qu71}]
The immersion $(f\times f)_G$ represents $P\big( f_*(1)\big)$; that is,
\[
P\big( f_*(1)\big)=\big((f\times f)_G\big)_*(1)\in h^{2dk}((X\times X)_G).
\]
\end{prop}

\section{Proof of Theorem 1}

We are now ready to prove Theorem 1. Recall that $f\co M^{n-dk}\imm X^n$ is assumed to be proper, self-transverse and $h^*$-oriented. We assume that $h^*$ has an external Steenrod operation $P$ of type $(G,d)$, and then the immersion $(f\times f)_G$ is $h^*$-oriented and represents $P\big( f_*(1)\big)$. The embedding of the self-intersection manifold
\[
\Sigma=\{(x,y)\mid x\neq y,f(x)=f(y)\}\hookrightarrow M\times M
\]
induces an embedding $\jmath\co\Sigma_G\hookrightarrow (M\times M)_G$ of homotopy orbit spaces. Consider the {\em equivariant double point immersion of $f$}
\[
\psi_G\co \Sigma_G\imm B\times X,\qquad [e,x,y]\mapsto \big([e],f(x)\big).
\]
The normal bundle of $\psi_G$ is isomorphic to the pullback $\jmath^*(\nu_f\times\nu_f)_G = \jmath^*\nu_{(f\times f)_G}$, and is oriented accordingly. Thus $\psi_G$ represents a class $(\psi_G)_*(1)\in h^{2dk}(B\times X)$.
\begin{lemma}\label{lem1}
The following diagram of proper immersions
\[
\xymatrix{
\Sigma_G\sqcup (B\times M) \ar[rr]^-{\jmath\,\sqcup\triangle_G} \ar[d]_{\psi_G\sqcup (1\times f)} & &(M\times M)_G \ar[d]^{(f\times f)_G} \\
B\times X \ar[rr]^-{\triangle_G} & & (X\times X)_G }
\]
is sub-cartesian, where $\sqcup$ denotes disjoint union. The excess bundle $E$ is zero over $\Sigma_G$ and isomorphic to $\gamma\otimes\nu_f$ over $B\times M$, where $\gamma$ denotes the real line bundle associated to the double cover $E\to B$.
\end{lemma}
\begin{proof}
The fact that $f$ is an immersion implies that $(f\times f)_G^{-1}(B\times X)=\Sigma_G\sqcup (B\times M)$, so that the square is a pullback. To prove it is sub-cartesian, it suffices to check the bundle condition on each component. On $\Sigma_G$ this is equivalent to self-transversality of $f$. On $B\times M$ it is equivalent to the condition
\[
d\triangle\left( TX_{f(x)}\right)\cap df(TM_x)\times df(TM_x) = (df\times df)\big(d\triangle (TM_x)\big)
\]
for each $x\in M$, which clearly holds. The excess bundle is zero over $\Sigma_G$ (again by self-transversality of $f$) and over $B\times M$ is determined by the exact sequence
\[
\xymatrix{
0 \ar[r] & B\times\nu_f \ar[r] & \triangle_G^*(\nu_f\times\nu_f)_G = (\nu_f\oplus\nu_f)_G \ar[r] & E \ar[r] & 0.
}
\]
Combined with the bundle isomorphism
\[
(\nu_f\oplus\nu_f)_G \to (B\times\nu_f)\oplus \gamma\otimes \nu_f,\qquad[e,v,w]\mapsto ([e],v+w)\oplus [e,1]\otimes (v-w),
\]
this shows that the excess bundle is as claimed.
\end{proof}

We now apply Proposition \ref{Ronga} to the sub-cartesian diagram of Lemma \ref{lem1} above, with $c=1\in h^0\big( (M\times M)_G\big)$ the unit element. This gives the equation
\[
\mathcal{P}\big(f_*(1)\big) = (\psi_G)_*(1) + (f\times 1)_*e(\gamma\otimes\nu_f) \in h^{2dk}(B\times X).
\]
To complete the proof it remains to show that $(\psi_G)_*(1)=2\times\Psi_*(1)$. The normal bundle $\nu_{\Psi}$ is isomorphic to the pullback $\kappa^*\nu_{(f\times f)_G}$, where the map
\[
\kappa\co \Sigma/G\to (M\times M)_G,\qquad \kappa[x,y] = [e,x,y]
\]
is independent of the choice of base point $e\in E$ up to homotopy. Note that $\psi_G$ factors as
\[
\xymatrix{
\Sigma_G:=E \times_G\Sigma \ar[r]^-{c} & B \times\Sigma/G \ar[r]^-{1\times\Psi} & B\times X,
}
\]
where the first map $c([e,x,y])=([e],[x,y])$ is a double cover. We claim that the orientation of $\psi_G$ as $\jmath^*\nu_{(f\times f)_G}$ agrees with its orientation as $c^*\nu_{1\times\Psi}\cong c^*p^*\nu_\Psi$ coming from the above factorisation, where $p$ denotes the projection $B\times\Sigma/G\to \Sigma/G$. This follows since $\jmath\simeq \kappa\circ p\circ c$. Hence by Lemma \ref{Um} we have
\begin{align*}
(\psi_G)_*(1)  & = (1\times\Psi)_*\circ c_*(1) \\
     & = (1\times\Psi)_*\circ c_*\circ c^*(1\times 1) \\
     &  = (1\times\Psi)_*(2\times 1)\\
     & = 2\times\Psi_*(1)
\end{align*}
as claimed.\hspace{\fill}$\Box$

\section{Ordinary Cohomology}

Of course, Steenrod operations were originally defined in ordinary cohomology. In this Section we interpret our formula for $h^*=H^*(-;\Z_2)$, the ordinary cohomology theory with coefficients in the ring $\Z_2$. This gives a new proof of a classical result of R.\ Thom \cite{Th50}.

Recall that the {\em Steenrod squares} are stable cohomology operations
\[
Sq^i\co H^*(-) \to H^{*+i}(-),\qquad i\geq 0
\]
satisfying various axioms. To construct the squares, one usually proceeds by first constructing an internal Steenrod operation of type $(G,1)$,
\[
\mathcal{P}\co H^*(-)\to H^{2*}(\R P^\infty\times -)
\]
To define the action of the squares on an element $\alpha\in H^k(X)$, note that the K\"{u}nneth Theorem gives an isomorphism $H^{2k}(\R P^\infty\times X)\cong \bigoplus_{i=0}^{2k} H^i(\R P^\infty)\otimes H^{2k-i}(X)$, and we have the formula
\[
\mathcal{P}(\alpha) = \sum_{i=0}^{k} \mu^{k-i}\otimes Sq^i(\alpha),\qquad\mu\in H^1(\R P^\infty)\mbox{ the generator}.
\]
The usual axioms then follow (see Steenrod and Epstein \cite{SE}). Note that for each $\ell\in \N$ we may restrict via the map $\iota_\ell\co \R P^\ell\times X\to \R P^\infty \times X$ to obtain a Steenrod operation
\[
\mathcal{P}_\ell\co H^*(-)\to H^{2*}(\R P^\ell\times -),\qquad \mathcal{P}_\ell = \iota_\ell^*\circ \mathcal{P}
\]
of type $(G,1)$ with respect to the principal $G$-bundle $S^\ell\to \R P^\ell$.

\begin{prop}
Let $f\co M^{n-k}\imm X^n$ be a proper self-transverse immersion representing $f_*(1)\in H^k(X)$. Then for all $\ell\in\N$,
\[
\mathcal{P}_\ell\big(f_*(1)\big)=(1\times f)_*
w_k(\gamma_\ell\otimes \nu_f)\in H^{2k}(\R P^\ell\times X),
\]
where $\gamma_\ell$ is the real line bundle associated to $S^\ell\to \R P^\ell$. Hence
\[
Sq^i\big(f_*(1)\big)=f_*w_i(\nu_f)\quad\textrm{for all }i.
\]
\end{prop}
\begin{proof}
The first statement follows immediately from Theorem 1 for the case $h^*=H^*(-;\Z_2)$, since in this theory $2=0$ and the Euler class is the top Stiefel-Whitney class.

To deduce the second statement from the first, note that for any $\alpha\in H^k(X)$ and $\ell >k$ we have
\[
\mathcal{P}_\ell( \alpha) = \sum_{i=0}^{k} \mu^{k-i}\otimes Sq^i(\alpha),
\]
where $\mu\in H^1(\R P^\ell)$ is the generator. On the other hand, the first statement gives
\begin{eqnarray*}
\mathcal{P}_\ell\big(f_*(1)\big) & = & (1\times f)_*
w_k(\gamma_\ell\otimes \nu_f)\\
    & = & (1\times f)_*\sum_{i=0}^k \mu^{k-i}\otimes w_i(\nu_f)\\
    & = & \sum_{i=0}^k \mu^{k-i}\otimes f_* w_i(\nu_f),
\end{eqnarray*}
by an easy calculation using the splitting principle and the fact that
$w_1(\eta\otimes\xi)=w_1(\eta)+w_1(\xi)$ for line bundles $\eta$ and
$\xi$ over the same base.
\end{proof}

\begin{rems}{\em When $h^*=\mathfrak{N}^*$ is unoriented cobordism, the same argument
shows that
\[
\mathcal{P}_\ell\big(f_*(1)\big)=(1\times f)_*
e(\gamma_\ell\otimes \nu_f)\in \mathfrak{N}^{2k}(\R P^\ell\times X).
\]
One may then use the results of \cite[Section 3]{Qu71} to derive a formula for $R^i(f_*(1))$ (where
the $R^i$ are the Steenrod operations in the theory $\mathfrak{N}^*$, defined in \cite[Section 15]{tD68}). The formula so obtained involves the coefficients of the formal group law in $\mathfrak{N}^*$, and the pushforwards  $f_*W_\alpha(\nu_f)$ of monomials in the Conner-Floyd Stiefel-Whitney classes of $\nu_f$
\cite{CF66}.
}\end{rems}

\end{document}